\documentclass{emsprocart}

\newcommand{\fA}{{\mathfrak{A}}}
\newcommand{\fS}{{\mathfrak{S}}}

\newcommand{\cE}{{\mathcal{E}}}
\newcommand{\cO}{{\mathcal{O}}}
\newcommand{\cP}{{\mathcal{P}}}

\newcommand{\bC}{{\mathrm{C}}}
\newcommand{\bG}{{\mathbf{G}}}
\newcommand{\bH}{{\mathbf{H}}}
\newcommand{\bL}{{\mathbf{L}}}
\newcommand{\bN}{{\mathrm{N}}}
\newcommand{\bT}{{\mathbf{T}}}

\newcommand{\FF}{{\mathbb{F}}}
\newcommand{\ZZ}{{\mathbb{Z}}}

\newcommand{\Aut}{\operatorname{Aut}}
\newcommand{\IBr}{\operatorname{IBr}}
\newcommand{\Irr}{\operatorname{Irr}}
\newcommand{\Ht}{\operatorname{ht}}
\newcommand{\mh}{\operatorname{mh}}
\newcommand{\GL}{\operatorname{GL}}
\newcommand{\PSL}{\operatorname{PSL}}
\newcommand{\SL}{\operatorname{SL}}
\newcommand{\Sp}{\operatorname{Sp}}
\newcommand{\GU}{\operatorname{GU}}
\newcommand{\SU}{\operatorname{SU}}
\newcommand{\ses}{{\operatorname{ss}}}
\newcommand{\RLG}{R_\bL^\bG}
\newcommand{\tw}[1]{{}^#1\!}

\let\eps=\epsilon
\let\la=\lambda

\contact[malle@mathematik.uni-kl.de]{Gunter Malle, Fachbereich Mathematik,
  TU Kaiserslautern, Postfach 3049, 67653 Kaiserslautern, Germany}

\numberwithin{equation}{section}

\newtheorem{theorem}{Theorem}[section]

\newtheorem{conjecture}[theorem]{Conjecture}

\theoremstyle{definition}

\newtheorem{problem}[theorem]{Problem}

\title[Local-global conjectures in representation theory]{Local-global conjectures in the\\ representation theory of finite groups}

\author[Gunter Malle]{Gunter Malle\thanks{The author gratefully acknowledges
financial support by the DFG priority program 1388.}}

\begin{document}

\begin{abstract}
We give a survey of recent developments in the investigation of the various
local-global conjectures for representations of finite groups.
\end{abstract}



\maketitle


\section{Introduction}

The \emph{ordinary} representation theory of finite groups over the complex
numbers was developed by Frobenius, Burnside, Schur and others around the
beginning of the 20th century, and the study of \emph{modular} representation
theory, dealing with
representations over fields of positive characteristic, was initiated by
Brauer around the mid of the 20th century. Still, amazingly enough, many
fundamental and easily formulated questions remain open in the ordinary as
well as in the modular representation theory of finite groups. For example,
in 1963 Brauer \cite{Br63} formulated a list of deep conjectures about
ordinary and modular representations of finite groups, most of which are yet
unsolved, and further important conjectures were subsequently put forward by
McKay, Alperin, Brou\'e and Dade. It is the purpose of this survey to expound
some of
the recent considerable advances on several of the major open conjectures in
this area. In this sense, this article can be seen as a continuation of the
1991 survey by Michler \cite{Mi91}. In particular, the opening question in
the introduction to that paper whether central conjectures in (modular)
representation theory might be provable as a consequence of the classification
of finite simple groups, has recently been given at least a partial positive
answer.
\par
We will concentrate here on so-called local-global conjectures which propose
to relate the representation theory of a finite group $G$ to that of its
\emph{local subgroups} for some prime $p$, that is, of subgroups $\bN_G(P)$
where $P\le G$ is a non-trivial $p$-subgroup. The charm of these conjectures,
as so often, lies in the stunning simplicity of their formulation as opposed
to their seeming intractability. More specifically, we will discuss the McKay
conjecture, its block-wise version known as Alperin--McKay conjecture,
Brauer's height zero conjecture and Dade's conjecture, all of which concern
character degrees of finite groups, as well as the Alperin weight conjecture,
which postulates a formula for the number of modular irreducible characters
in terms of local data. A possible structural explanation of some instances of
these conjectures is offered by Brou\'e's abelian defect group conjecture.
All of these conjectures and observations point towards some hidden theory
explaining these phenomena, but we are unable to find it yet.
\par
The approach of using local data to obtain global information on a finite group
had already proved very successful in the classification of finite simple
groups. Now, conversely, the classification seems to provide a way for proving
local-global conjectures in representation theory. The basic idea is to reduce
the conjectures to possibly more complicated statements about finite simple
groups, which we then hope to verify by using our detailed knowledge on these
groups. This approach has already proved successful in two
important cases, see Theorems~\ref{thm:MS} and~\ref{thm:H0}.
\par
Not unexpectedly, the attempt to apply the classification has on the one hand
led to the development of new, purely representation theoretic notions,
tools and results, and on the other hand it has made apparent that our
knowledge even of the ordinary representation theory of the finite simple
groups is far from sufficient for many purposes. In this way, the reduction
approach to the local-global conjectures has already spawned powerful new
methods, interesting questions and challenging research topics even outside
its immediate range.
\vskip 1pc
\noindent
{\bf Acknowledgement:} The author thanks Olivier Dudas, Eugenio Giannelli,
Radha Kessar, Caroline Lassueur, Gabriel Navarro and Britta Sp\"ath for
various helpful comments on an earlier version.

\section{The McKay Conjecture, the Alperin--McKay Conjecture and refinements}

The McKay conjecture \cite{Mc72} is the easiest local-global conjecture to
state. It could already have been formulated by Frobenius or Burnside, but
was first noticed only in 1972. It is also the origin, together with Alperin's
Weight Conjecture~\ref{conj:BAW}, of the more general Dade
Conjecture~\ref{conj:Dade} as well as of Brou\'e's Conjecture~\ref{conj:Br}.

\subsection{Characters of $p'$-degree}
For a finite group $G$ and a prime $p$ we denote by
$$\Irr_{p'}(G):=\{\chi\in\Irr(G)\mid p\text{ does not divide }\chi(1)\}$$
the set of irreducible complex characters of $G$ of degree prime to $p$.

\begin{conjecture}[McKay (1972)]   \label{conj:McK}
 Let $G$ be a finite group and $p$ be a prime. Then
 $$|\Irr_{p'}(G)|=|\Irr_{p'}(\bN_G(P))| \, ,$$
 where $P$ is a Sylow $p$-subgroup of $G$ and $\bN_G(P)$ its normaliser.
\end{conjecture}

That is to say, certain fundamental information on the representation theory
of $G$ is encoded in local subgroups of $G$, namely in the Sylow normalisers.

In fact, McKay \cite{Mc72} made his conjecture only for $G$ a simple group 
and for the prime $p=2$. It was Isaacs, in his landmark paper \cite{Is73},
who proved the conjecture for all groups of odd order and any prime $p$.
Soon afterwards, Alperin \cite{Al76} refined and extended the statement of
Conjecture~\ref{conj:McK} to include Brauer blocks, now known as the
Alperin--McKay conjecture. To formulate it let us fix a $p$-modular system
$(K,\cO,k)$, where $\cO$ is a discrete valuation ring with field of fractions
$K$ of characteristic~0 and with finite residue field $k$ of characteristic~$p$,
large enough for the given finite group $G$. Then the group ring $\cO G$
decomposes as a direct sum of minimal 2-sided ideals, the \emph{$p$-blocks} of
$G$, and every irreducible character of $G$ is non-zero on exactly one of these
blocks. This induces a partition $\Irr(G)=\coprod_B \Irr(B)$ of the irreducible
characters of $G$, where $B$ runs
over the $p$-blocks of $G$. To each block $B$ is attached a $p$-subgroup
$D\le G$, uniquely determined up to conjugacy, a so-called \emph{defect group}
of $B$. For a block $B$ with defect group $D$ we then write
$$\Irr_0(B):=\{\chi\in\Irr(B)\mid \Ht(\chi)=0\}$$
for the set of \emph{height zero characters} in $B$; here the \emph{height}
$\Ht(\chi)$ of the irreducible character $\chi$ is defined by the formula
$\chi(1)_p|D|_p=p^{\Ht(\chi)}|G|_p$. Thus $\Irr_0(B)=\Irr_{p'}(B)$ if $D$ is
a Sylow $p$-subgroup of $G$.
Brauer has shown how to construct a $p$-block $b$ of $\bN_G(D)$, closely
related to $B$, called the \emph{Brauer correspondent} of $B$. We then also
say that $B=b^G$ is the Brauer induced block from $b$.

\begin{conjecture}[Alperin (1976)]   \label{conj:AM}
 Let $G$ be a finite group, $p$ be a prime and $B$ a $p$-block of $G$ with
 defect group $D$. Then
 $$|\Irr_0(B)|=|\Irr_0(b)| \, ,$$
 where $b$ is the Brauer correspondent of $B$ in $\bN_G(D)$.
\end{conjecture}

Clearly, by summing over all blocks of maximal defect, that is, blocks whose
defect groups are Sylow $p$-subgroups of $G$, the Alperin--McKay
Conjecture~\ref{conj:AM} implies the McKay Conjecture~\ref{conj:McK}.

Soon after its formulation the Alperin--McKay Conjecture~\ref{conj:AM} was
proved  for $p$-solvable groups by Okuyama and Wajima \cite{OW80} and
independently by Dade \cite{D80}. It has also been verified for symmetric
groups $\fS_n$ and alternating groups $\fA_n$ by Olsson \cite{Ol76}, and for
their covering groups and for the general linear groups $\GL_n(q)$ by Michler
and Olsson \cite{MO83,MO90}.

Subsequently, several refinements of this conjecture were proposed. The first
one by Isaacs and Navarro \cite{IN02} predicts additional congruences; here
$n_{p'}$ denotes the part prime to $p$ of an integer~$n$:

\begin{conjecture}[Isaacs--Navarro (2002)]   \label{conj:IN}
 In the situation of Conjecture~\ref{conj:AM} there exists a bijection
 $\Omega:\Irr_0(B)\rightarrow\Irr_0(b)$ and a collection of signs
 $(\eps_\chi|\chi\in\Irr_0(B))$ such that
 $$\Omega(\chi)(1)_{p'}\equiv\eps_\chi\chi(1)_{p'}\pmod{p}.$$
\end{conjecture}

(Note that this is a true refinement of Conjecture~\ref{conj:AM} whenever
$p\ge5$.)
This has been shown to hold for example for $\fS_n$, $\fA_n$ and their double
covers by Fong \cite{F03}, Nath \cite{N09} and Gramain \cite{Gr11}
respectively. Two further
refinements on the properties of the required bijection concerning the action
of those Galois automorphisms fixing a prime ideal above $p$ were put forward
in the same paper \cite{IN02}, and by Navarro \cite{Na04} respectively.
Yet another refinement due to Turull \cite{Tu07} includes $p$-adic fields and
Schur indices. 

\subsection{A reduction theorem}
While Conjecture~\ref{conj:McK} was subsequently checked for several
further families of finite groups, the first significant breakthrough in the
case of general groups was achieved by Isaacs, Navarro and the author
\cite{IMN} in 2007 where they reduced the McKay conjecture to a question on
simple groups:

\begin{theorem}[Isaacs--Malle--Navarro (2007)]   \label{thm:IMN}
 The McKay Conjecture~\ref{conj:McK} holds for all finite groups at the prime
 $p$, if all finite \emph{non-abelian simple} groups satisfy the so-called
 \emph{inductive McKay condition} at the prime $p$.
\end{theorem}

This inductive condition for a simple group $S$ is stronger than just the
validity of McKay's conjecture for $S$, and in particular also involves the
covering groups and the automorphism group of the simple group in question:
If $S$ is simple, and $G$ is its universal covering group (the largest perfect
central extension of $S$), then the \emph{inductive McKay condition} on $S$ is
satisfied, if for some proper $\Aut(G)_P$-invariant subgroup $M<G$ containing
the normaliser $\bN_G(P)$ of a Sylow $p$-subgroup $P$ of $G$
\begin{enumerate}
\item[(1)] there exists an $\Aut(G)_M$-equivariant bijection
  $\Omega:\Irr_{p'}(G)\rightarrow\Irr_{p'}(M)$, respecting central characters,
\item[(2)] such that the extendibility obstructions of $\chi$ and
  $\Omega(\chi)$ to their respective inertia groups in $G\rtimes\Aut(G)$,
  considered as 2-cocycles, coincide for all $\chi\in\Irr_{p'}(G)$.
\end{enumerate}
Here, for $X\le G$, $\Aut(G)_X$ denotes the stabiliser of $X$ in $\Aut(G)$.
Note that due to the inductive nature of the reduction argument we need not
descend all the way to $\bN_G(P)$, but only to some intermediary subgroup~$M$
of our choice. As will be seen below this is very useful in the case of finite
groups of Lie type. In fact, the condition stated in \cite[\S10]{IMN} (where
this notion is called \emph{being good for the prime $p$}) even allows for a
slightly bigger group to be considered in place of $G$, which is particularly
useful in dealing with the finite groups of Lie type.
\par
The inductive McKay condition has been shown for the alternating and sporadic
groups by the author \cite{Ma08}, and for groups of Lie type in their defining
characteristic by Sp\"ath \cite{Sp12}, extending work of Brunat \cite{Br09}
and building on a result of Maslowski \cite{Ms10}. Thus only the simple groups
of Lie type at primes $p$ different from their defining characteristic remain
to be considered. These are best studied as finite reductive groups.

\subsection{The inductive condition for groups of Lie type}
If $G$ is the universal covering group of a simple group of Lie type, then up
to finitely many known exceptions
(see e.g.~\cite[Tab.~24.3]{MT}) there exists a simple linear algebraic group
$\bG$ of simply connected type over the algebraic closure of a finite field
$\FF_q$ and a Steinberg endomorphism $F:\bG\rightarrow\bG$ such that
$G=\bG^F$ is the finite group of fixed points in $\bG$ under $F$, a finite
reductive group. Lusztig has obtained a parametrisation of the irreducible
complex characters of the groups $\bG^F$ and in particular has determined their
degrees. To describe this, let's assume for simplicity that $F$ is the
Frobenius map with respect to some $\FF_q$-structure of $\bG$. Let $\bG^*$ be
a \emph{dual group to} $\bG$ (with root datum dual to the one of $\bG$)
and with compatible Frobenius map on $\bG^*$ also denoted by $F$. Then Lusztig
\cite{LuB} has constructed a partition
$$\Irr(G)=\coprod_{s\in G_\ses^*/\sim}\cE(G,s)$$
of the irreducible characters of $G$ into \emph{Lusztig series} $\cE(G,s)$
parametrised
by semisimple elements $s\in G^*:=\bG^{*F}$ up to $G^*$-conjugacy. Further for
any semisimple element $s\in G^*$ he obtained a \emph{Jordan decomposition}
$$\Psi_s:\cE(G,s)\buildrel{1 - 1}\over\longrightarrow\cE(\bC_{G^*}(s),1)$$
relating the Lusztig series $\cE(G,s)$ to the so-called \emph{unipotent
characters} of the (possibly disconnected) group $\bC_{G^*}(s)$, such that
the degrees satisfy
\begin{align}   
\chi(1)& = |G^*:\bC_{G^*}(s)|_{q'}\,\Psi_s(\chi)(1)\qquad
  \text{ for all }\chi\in\cE(G,s).\label{eq:JD}
\end{align}
The unipotent characters of finite reductive groups have been classified by
Lusztig \cite{LuB} and he has given combinatorial formulas for their degrees.
It is thus in principle possible to determine the irreducible characters of
$G$ of $p'$-degree. For example, if $p$ is a prime not dividing $q$,
Equation~(\ref{eq:JD}) shows that $\chi\in\cE(G,s)$ lies in $\Irr_{p'}(G)$ if
and only if $s$ centralises a Sylow $p$-subgroup of $G^*$ and the Jordan
correspondent $\Psi_s(\chi)$ lies in $\Irr_{p'}(\bC_{G^*}(s))$, thus yielding
a reduction to unipotent characters.
\par
The proper tool for discussing unipotent characters is provided by
\emph{$d$-Harish-Chandra theory}, introduced by Brou\'e--Malle--Michel
\cite{BMM} and further developed by the author \cite{MaU,Ma00,Ma07}. For this,
let
$$d=d_p(q):=\text{multiplicative order of $q$ }
  \begin{cases} \text{modulo $p$}& \text{ if $p$ is odd},\\
               \text{modulo $4$}& \text{ if $p=2$}.\end{cases}$$
In \cite{Ma07} we give a parametrisation of $\Irr_{p'}(G)$ in terms of
combinatorial data related to the \emph{relative Weyl group}
$\bN_G(\bT_d)/\bC_G(\bT_d)$ of a Sylow $d$-torus $\bT_d$ of $G$. This is
always a finite complex reflection group. Here an $F$-stable torus $\bT\le\bG$
is called a \emph{$d$-torus} if it splits over $\FF_{q^d}$ and no $F$-stable
subtorus of $\bT$ splits over any smaller field. A $d$-torus of maximal
possible dimension in $\bG$ is called a \emph{Sylow $d$-torus}. Such Sylow
$d$-tori are unique up to $G$-conjugacy~\cite{BM92}.

On the other hand, the following result \cite[Thms.~5.14 and~5.19]{Ma07} shows
that in most cases we may choose $M:=\bN_G(\bT_d)$ as the intermediary
subgroup occurring in the inductive McKay condition:

\begin{theorem}[Malle (2007)]  \label{thm:normSyl}
 Let $\bG$ be simple, defined over $\FF_q$ with corresponding Frobenius map
 $F:\bG\rightarrow\bG$ and let $G:=\bG^F$. Let $p{\not|}q$ be a prime divisor
 of $|G|$, and set $d=d_p(q)$. Then the normaliser $\bN_G(\bT_d)$ of a Sylow
 $d$-torus $\bT_d$ of $G$ contains a Sylow $p$-subgroup of $G$ unless one of
 the following holds:
 \begin{itemize}
 \item[(a)] $p=3$, and $G=\SL_3(q)$ with $q\equiv4,7\pmod9$, $G=\SU_3(q)$
  with $q\equiv2,5\pmod9$, or $G=G_2(q)$ with $q\equiv2,4,5,7\pmod9$; or
 \item[(b)] $p=2$, and $G=\Sp_{2n}(q)$ with $n\ge1$ and $q\equiv3,5\pmod8$.
 \end{itemize}
\end{theorem}

In particular, with this choice $M$ only depends on $d$, but not on the precise
structure of a Sylow $p$-subgroup or the Sylow $p$-normaliser, which makes a
uniform
argument feasible. The four exceptional series in Theorem~\ref{thm:normSyl}
can be dealt with separately (see \cite{Ma08b}). For example, part~(b) includes
the case that $G=\SL_2(q)$ where $q\equiv3,5\pmod8$ and the Sylow 2-normaliser
is isomorphic to $\SL_2(3)$, while torus normalisers are dihedral groups.
For the general case, in a delicate
Clifford theoretic analysis Sp\"ath \cite{Sp09,Sp10} has shown that
$\Irr_{p'}(\bN_G(\bT_d))$ can be parametrised by the same combinatorial
objects as for $\Irr_{p'}(G)$, thus completing the proof of:

\begin{theorem}[Malle (2007) and Sp\"ath (2010)]
 Let $\bG$ be simple, defined over $\FF_q$ with corresponding Frobenius map
 $F:\bG\rightarrow\bG$ and let $G:=\bG^F$. Let $p{\not|}q$ be a prime divisor
 of $|G|$, $d=d_p(q)$, and assume that we are not in one of the exceptions (a)
 or~(b) of Theorem~\ref{thm:normSyl}. Then there is a bijection
 $$\Omega:\Irr_{p'}(G)\rightarrow\Irr_{p'}(\bN_G(\bT_d))\ \text{with }
   \Omega(\chi)(1)\equiv\pm\chi(1)\!\!\!\!\pmod p\text{ for }
   \chi\in \Irr_{p'}(G).$$
\end{theorem}

So in particular we obtain degree congruences as predicted by
Conjecture~\ref{conj:IN}.

The equivariance and cohomology properties of such a bijection
$\Omega$ required by the inductive McKay condition have at present
been shown by Cabanes--Sp\"ath \cite{CS13,CS15b} and the author \cite{Ma08b}
for all series of groups of Lie type except types $B_n$, $C_n$, $D_n$,
$\tw2D_n$, $E_6$, $\tw2E_6$
and $E_7$. The most difficult and complicated part was certainly the proof by
Cabanes and Sp\"ath \cite{CS15b} that the linear and unitary groups do
satisfy the inductive McKay condition. It relies on a powerful criterion
of Sp\"ath which allows one to descend a bijection for the much easier case
of $\GL_n(q)$, for example, to its quasi-simple subgroup $\SL_n(q)$ if the
inertia groups of $p'$-characters of $\SL_n(q)$ and of the intermediary
subgroup $M$ have a certain semidirect product decomposition, see
\cite[Thm.~2.12]{Sp12} for details. This criterion is shown to hold for linear
and unitary groups using Kawanaka's generalised Gelfand Graev characters.
\par
The treatment of the remaining seven series of groups seems to require further
knowledge on their ordinary representation theory, not immediate from Lusztig's
results. More precisely, a solution will need to solve the following:

\begin{problem}
 For $G$ quasi-simple of Lie type, determine the action of $\Aut(G)$ on
 $\Irr(G)$. 
\end{problem}

More precisely, it is not known in general how outer automorphisms act on
irreducible complex characters of $G$ lying in series $\cE(G,s)$ with
$\bC_{\bG^*}(s)$ not connected. In particular, the ordinary character degrees
of extensions of $G$ by outer automorphisms are unknown.

The most recent and most far-reaching result in this area has been obtained by
Sp\"ath and the author \cite{MS15}, showing that McKay's original question has
an affirmative answer:

\begin{theorem}[Malle--Sp\"ath (2015)]   \label{thm:MS}
 The McKay conjecture holds for all finite groups at the prime $p=2$.
\end{theorem}

For the proof we show that the groups in the remaining seven families also
satisfy the inductive McKay condition at the prime $2$ and then apply
Theorem~\ref{thm:IMN}. This relies on an equivariant extension of the
Howlett--Lehrer theory of endomorphism algebras of induced cuspidal modules
for finite groups with a BN-pair, and on special properties of the prime~2
as a divisor of character degrees of groups of Lie type. Namely, except for the
characters of degree $(q-1)/2$ of $\SL_2(q)$ for $q\equiv3\pmod4$, non-linear
cuspidal characters are always of even degree. The latter statement fails
drastically for odd primes. An immediate extension to other primes thus seems
very hard.

The result of Theorem~\ref{thm:MS} shows that the approach to the local-global
conjectures via the reduction to finite simple groups is indeed successful. 

\subsection{The block-wise reduction}
The strategy and proof of Theorem~\ref{thm:IMN} have become the blueprint
for all later reductions of other local-global conjectures. So Sp\"ath
\cite{Sp13a} saw how this reduction could be (simplified and then) extended
to the block-wise setting:

\begin{theorem}[Sp\"ath (2013)]
 The Alperin--McKay Conjecture~\ref{conj:AM} holds for all finite groups at
 the prime $p$, if all finite \emph{non-abelian simple} groups satisfy the
 so-called \emph{inductive Alperin--McKay condition} at $p$.
\end{theorem}

In fact, her reduction also applies to the more precise Isaacs--Navarro
Conjecture~\ref{conj:IN} involving degree congruences.

The \emph{inductive Alperin--McKay condition} on a simple group $S$ is quite
similar to the inductive McKay condition as outlined above: Let $G$ denote the
universal covering group of $S$. Then for each isomorphism class of defect
group $D\le G$ we need a subgroup $\bN_G(D)\le M_D\le G$, proper unless $D$ is
central, such that for each block $B$ with defect group $D$ and Brauer
corresponding block $b$ of $M_D$ there exists an $\Aut(G)_b$-equivariant
bijection $\Omega:\Irr_0(B)\rightarrow\Irr_0(b)$, respecting central
characters and having further rather technical properties phrased in terms of
projective representations of automorphism groups of $G$, see
\cite[Def.~7.2]{Sp13a} for details, as well as the article by Sp\"ath in this
volume \cite{Sp16}. Here again $\Aut(G)_b$ denotes the
stabiliser of the block $b$ in $\Aut(G)$. This condition has been verified by
Koshitani and Sp\"ath \cite{KS15a,KS15b} for all blocks with cyclic defect
groups, as well as for groups of Lie type in their defining characteristic and
alternating groups at odd primes by Sp\"ath \cite{Sp13a}, while Denoncin
\cite{De14} proved it for alternating groups at $p=2$. Cabanes and Sp\"ath
\cite{CS15a} show it for blocks of $\SU_n(q)$ and $\SL_n(q)$ of maximal
defect. For the sporadic groups see the
website by Breuer \cite{BrWeb}.

In this context we mention the following open question:

\begin{problem}
 Find a reduction for the Alperin--McKay conjecture including the action of
 certain Galois automorphisms as predicted by Isaacs and Navarro
 \cite{IN02,Na04}.
\end{problem}
 
A recent result of Ladisch \cite{La15} can be seen as a first step towards such
a reduction. This might also give a hint for even more natural bijections in
the verification of the inductive conditions for groups of Lie type.

\section{Brauer's Height Zero Conjecture}

The Alperin--McKay Conjecture~\ref{conj:AM} predicts the number of characters
of height zero by local data. When are these all the irreducible characters
in a given block?

\subsection{Characters of height zero}
An answer is postulated in Brauer's Height Zero Conjecture \cite{Br55}
from~1955:

\begin{conjecture}[Brauer (1955)]   \label{conj:H0}
 Let $B$ be a $p$-block of a finite group with defect group $D$. Then
 all irreducible characters in $B$ have height zero if and only if $D$ is
 abelian.
\end{conjecture}

A positive solution would provide, for example, an extremely
simple method to detect from a group's ordinary character table whether its
Sylow $p$-subgroups are abelian: indeed, the Sylow $p$-subgroups are defect
groups of the principal block, and the characters (degrees) in the latter can
be read off from the character table.

The $p$-solvable case of Conjecture~\ref{conj:H0} is an impressive theorem
by Gluck and Wolf \cite{GW84}. All further substantial progress on this
question was made using the classification of finite simple groups. The most
far reaching general result so far concerns 2-blocks whose defect groups are
Sylow 2-subgroups \cite{NT12}:

\begin{theorem}[Navarro--Tiep (2012)]   \label{thm:maxdef}
 Let $B$ be a 2-block of a finite group of maximal defect. Then Brauer's
 Height Zero Conjecture~\ref{conj:H0} holds for $B$.
\end{theorem}

In particular, the above criterion for detection of abelian Sylow $p$-subgroups
from the ordinary character table holds when $p=2$.

The proof of Theorem~\ref{thm:maxdef} relies on Walter's determination of
finite groups with abelian Sylow 2-subgroups as well as on Lusztig's
previously described classification of irreducible characters of finite
reductive groups. 

\subsection{The ``if'' direction}
For the case of arbitrary blocks and primes, Berger and Kn\"orr \cite{BK}
derived the following optimal reduction to the same statement for blocks of
quasi-simple groups (recall that a finite group $G$ is \emph{quasi-simple}
if $G$ is perfect and $G/Z(G)$ is simple):

\begin{theorem}[Berger--Kn\"orr (1988)]   \label{thm:BK}
 The ``if''-direction  of Brauer's Height Zero Conjecture~\ref{conj:H0} holds
 for the $p$-blocks of all finite groups, if it holds for the $p$-blocks of all
 quasi-simple groups.
\end{theorem}

First significant steps in the verification of the assumption of this reduction
theorem were subsequently obtained by Olsson \cite{Ol90} for the covering
groups of alternating groups. The case of groups of Lie type in their
defining characteristic is easy for this questions, as defect groups are
either Sylow $p$-subgroups or trivial, and Sylow $p$-subgroups are non-abelian
unless we are in the case of $\PSL_2(q)$. For non-defining characteristic,
Blau and Ellers \cite{BE} obtained the following important result:

\begin{theorem}[Blau--Ellers (1999)]   \label{thm:BE}
 Brauer's Height Zero Conjecture~\ref{conj:H0} holds for all blocks of
 quasi-simple central factor groups of $\SL_n(q)$ and $\SU_n(q)$.
\end{theorem}

\subsection{Blocks of groups of Lie type}
The case of the other quasi-simple groups of Lie type could only be settled
after having obtained a full parametrisation of their $p$-blocks. This
classification is very closely related to Lusztig induction and
can again be most elegantly phrased in terms of $d$-Harish-Chandra theory.
It was achieved over a period of over 30 years by work of many authors.
As before let $\bG$ be a connected reductive algebraic group defined over
$\FF_q$ with corresponding
Frobenius endomorphism $F:\bG\rightarrow\bG$, and let $\bG^*$ be dual to $\bG$.
The first general reduction step was given by Brou\'e and Michel \cite{BM89}
who showed a remarkable compatibility between Brauer blocks and Lusztig series:
for any semisimple $p'$-element $s\in \bG^{*F}$ the set
$\cE_p(\bG^F,s):=\coprod_t \cE(\bG^F,st)$ is a union of $p$-blocks, where
$t$ runs over $p$-elements in the centraliser $\bC_{\bG^*}(s)^F$. All further
progress is linked to Lusztig induction. For an $F$-stable Levi subgroup
$\bL\le\bG$, using $\ell$-adic cohomology of suitable varieties attached to
$\bL$ and $\bG$, Lusztig has defined an induction map
$$\RLG:\ZZ\Irr(\bL^F)\rightarrow\ZZ\Irr(\bG^F).$$
Proving a conjecture of Brou\'e, Bonnaf\'e and Rouquier \cite{BR03}
showed that most of the series $\cE_p(\bG^F,s)$ ``come from below'' (see
also the recent extension of this result by Bonnaf\'e , Dat and Rouquier
\cite[Thm.~7.7]{BDR}):

\begin{theorem}[Bonnaf\'e--Rouquier (2003)]   \label{thm:BR}
 Let $s\in\bG^{*F}$ be a semisimple $p'$-element, and let $\bL\le\bG$ be an
 $F$-stable Levi subgroup such that $\bC_{\bG^*}(s)\le\bL^*$. Then $\RLG$
 lifts to Morita equivalences between the blocks in $\cE_p(\bL^F,s)$ and in
 $\cE_p(\bG^F,s)$.
\end{theorem}

This reduces the determination of blocks to the so-called \emph{quasi-isolated}
situation, that is to series $\cE_p(\bG^F,s)$ where $\bC_{\bG^*}(s)$ is not
contained in any proper $F$-stable Levi subgroup of $\bG^*$. Here crucial
steps were provided by Fong--Srinivasan \cite{FS86} for groups of
classical type, Brou\'e--Malle--Michel \cite{BMM} for unipotent blocks and
large primes, Cabanes--Enguehard \cite{CE99} for general blocks and primes
$p\ge5$, Enguehard \cite{E00} for unipotent blocks of exceptional type groups,
and Kessar--Malle \cite{KM13,KM15a} for the remaining quasi-isolated cases.
To describe the result, let
$$\cE(\bG^F,p'):=\{\chi\in\cE(\bG^F,s)\mid
   s\in\bG_\ses^{*F}\text{ of $p'$-order}\},$$
the set of irreducible characters lying in Lusztig series labelled by
$p'$-elements. Then $\RLG$ restricts to a map
$\ZZ\cE(\bL^F,p')\rightarrow\ZZ\cE(\bG^F,p')$. Levi subgroups of the form
$\bC_\bG(\bT)$, where $\bT$ is a $d$-torus of $\bG$ are called \emph{$d$-split},
and $\chi\in\Irr(\bG^F)$ is called \emph{$d$-cuspidal} if it does not occur as
a constituent of $\RLG(\la)$ for any proper $d$-split Levi subgroup $\bL<\bG$
and any $\la\in\Irr(\bL^F)$. More generally $\chi\in\cE(\bG^F,s)$ is called
\emph{$d$-Jordan cuspidal} if its Jordan correspondent
$\Psi_s(\chi)\in\cE(C_{\bG^*}(s)^F,1)$ is $d$-cuspidal. With this, the
classification of $p$-blocks (in the smoothest case) can be formulated as
follows in terms of Lusztig induction:

\begin{theorem}   \label{thm:p-blocks}
 Let $\bH$ be a simple algebraic group of simply connected type defined over
 $\FF_q$ with corresponding Frobenius endomorphism $F:\bH\rightarrow\bH$. Let
 $\bG\le \bH$ be an $F$-stable Levi subgroup. Let $p{\not|}q$ be a prime
 and set $d=d_p(q)$.
 \begin{enumerate}
  \item[\rm(a)] For any $d$-split Levi subgroup $\bL\le\bG$ and any
   $d$-Jordan-cuspidal character $\la\in\cE(\bL^F,p')$, there exists a unique
   $p$-block $b(\bL,\la)$ of $\bG^F$ such that all irreducible constituents
   of $\RLG(\la)$ lie in $b(\bL,\la)$.
  \item[\rm(b)] The induced map $(\bL,\la)\mapsto b(\bL,\la)$ on
   $G$-conjugacy classes of pairs as in~(a) is bijective if $p\geq 3$ is good
   for $\bG$, and if moreover $p\ne 3$ if $\bG^F$ has a factor $\tw3D_4(q)$.
 \end{enumerate}
\end{theorem}

A statement in full generality can be found in  \cite[Thm.~A]{KM15a}. Kessar
and the author \cite{KM13} used this classification to complete the proof of
the ``if'' direction of Brauer's Height Zero Conjecture~\ref{conj:H0}, relying
on the Berger--Kn\"orr reduction (Theorem~\ref{thm:BK}) and on the Blau--Ellers
result (Theorem~\ref{thm:BE}), thus offering further proof of the viability
of the reduction approach to local-global conjectures:

\begin{theorem}[Kessar--Malle (2013)]   \label{thm:H0}
 Let $B$ be a $p$-block of a finite group. If $B$ has abelian defect groups,
 then all irreducible characters in $B$ have height zero.
\end{theorem}

As an important step in the proof we show that the Bonnaf\'e--Rouquier Morita
equivalences in Theorem~\ref{thm:BR} preserve abelianity of defect groups (this
has now been reproved more conceptually in \cite{BDR}).

Navarro, Solomon and Tiep \cite{NST15} use Theorem~\ref{thm:H0} to derive an
effective criterion to decide the abelianity of Sylow subgroups from the
character table.

\subsection{The ``only if'' direction}
A crucial ingredient of Navarro and Tiep's proof of Theorem~\ref{thm:maxdef}
was a theorem of Gluck and Wolf for the prime~2. The missing odd-$p$ analogue
of this seemed to constitute a major obstacle towards establishing the
remaining, ``only if'' direction of the Height Zero Conjecture. Using the
classification of finite simple groups Navarro and Tiep \cite{NT13} have now
obtained a proof of this result:

\begin{theorem}[Navarro--Tiep (2013)]
 Let $N\unlhd G$ be finite groups, $p$ a prime, and $\theta\in\Irr(N)$ a
 $G$-invariant character. If $\chi(1)/\theta(1)$ is prime to $p$ for all
 $\chi\in\Irr(G)$ lying above $\theta$ then $G/N$ has abelian Sylow
 $p$-subgroups.
\end{theorem}

Building on this, Navarro and Sp\"ath \cite{NS14} succeeded in proving the
following reduction theorem for this direction of the conjecture: 

\begin{theorem}[Navarro--Sp\"ath (2014)]   \label{thm:NS}
 The ``only if''-direction of Brauer's Height Zero Conjecture~\ref{conj:H0}
 holds for all finite groups at the prime $p$, if
 \begin{enumerate}
  \item[(1)] it holds for all $p$-blocks of all quasi-simple groups, and
  \item[(2)] all simple groups satisfy the inductive Alperin--McKay condition
   at $p$.
 \end{enumerate}
\end{theorem}

For their proof, they introduce and study the new notion of central block
isomorphic character triples.

The first assumption of Theorem~\ref{thm:NS} was recently shown to hold
\cite{KM15b}, again building on the classification of blocks of finite
reductive groups described before:

\begin{theorem}[Kessar--Malle (2015)]
 The ``only if''-direction of Brauer's Height Zero Conjecture~\ref{conj:H0}
 holds for all $p$-blocks of all quasi-simple groups.
\end{theorem}

Thus, Brauer's height zero conjecture will follow once the inductive
Alperin--McKay condition has been verified for all simple groups. This
again underlines the central importance of the Alperin--McKay
Conjecture~\ref{conj:AM} in the representation theory of finite groups.

\subsection{Characters of positive height}
Conjecture~\ref{conj:H0} only considers characters of height zero. It is
natural to ask what can be said about the heights of other characters
in a given block. There are two conjectural answers to this question. To state
the first one, for $B$ a $p$-block we define
$$\mh(B):=\min\{\Ht(\chi)\mid \chi\in\Irr(B)\setminus\Irr_0(B)\},$$
the minimal positive height of a character in $\Irr(B)$, and we formally set
$\mh(B)=\infty$ if all characters in $B$ are of height~0. Eaton and Moret\'o
\cite{EM14} have put forward the following:

\begin{conjecture}[Eaton--Moret\'o (2014)]   \label{conj:EM}
 Let $B$ be a $p$-block of a finite group with defect group $D$. Then
 $\mh(B)=\mh(D)$.
\end{conjecture}

The case when $\mh(B)=\infty$ is Brauer's height zero conjecture, since
clearly all characters of the defect group $D$ are of height zero if and only
$D$ is abelian.

Eaton and Moret\'o \cite{EM14} proved their conjecture for all blocks
of symmetric and sporadic groups, and for $\GL_n(q)$
for the defining prime. They also showed that for $p$-solvable groups we
always have $\mh(D)\le\mh(B)$, and that this inequality is true for all groups
if Dade's projective conjecture (see Section~\ref{sec:Dade}) holds. Brunat
and the author \cite{BM15} then checked that the Eaton--Moret\'o
Conjecture~\ref{conj:EM} holds for all principal blocks of quasi-simple groups,
for all $p$-blocks of quasi-simple groups of Lie type in characteristic~$p$,
all unipotent blocks of quasi-simple exceptional groups of Lie type, and all
$p$-blocks of covering groups of an alternating or symmetric group. No
reduction of this conjecture to simple groups is known, though.

A different approach to characters of positive height is given by Dade's
Conjecture, which we review in Section~\ref{sec:Dade} below.

\section{The Alperin Weight Conjecture}   \label{sec:AWC}

While the McKay Conjecture counts characters of $p'$-degree, the Alperin
Weight Conjecture concerns characters whose degree has maximal possible
$p$-part, the so-called defect zero characters.

\subsection{Weights and chains} An irreducible  character $\chi$ of a finite
group $G$ has \emph{defect zero} if $\chi(1)_p=|G|_p$. A \emph{$p$-weight}
of $G$ is a pair $(Q,\psi)$ where $Q\le G$ is a \emph{radical $p$-subgroup},
that is, $Q=O_p(\bN_G(Q))$, and $\psi\in\Irr(\bN_G(Q)/Q)$ is a defect zero
character. If $\psi$ lies in the
block $b$ of $\bN_G(Q)$, then the weight $(Q,\psi)$ is said to \emph{belong to
the block $b^G$} of $G$. Alperin's original formulation of the weight
conjecture \cite{Al87} now proposes to count the $p$-modular irreducible
Brauer characters $\IBr(B)$ in a $p$-block $B$ in terms of weights:

\begin{conjecture}[Alperin (1986)]   \label{conj:BAW}
 Let $G$ be a finite group, $p$ be a prime and $B$ a $p$-block of $G$. Then
 $$|\IBr(B)|=|\{[Q,\psi]\mid (Q,\psi)\text{ a $p$-weight belonging to }B\}|,$$
 where $[Q,\psi]$ denotes the $G$-conjugacy class of the $p$-weight
 $(Q,\psi)$.
\end{conjecture}

The name ``weights'' was apparently chosen since for groups of Lie type in
defining characteristic the irreducible Brauer characters are indeed labelled
by (restricted) weights of the corresponding linear algebraic groups.

Alperin \cite{Al87} notes the following nice consequence of his conjecture:

\begin{theorem}[Alperin]
 Assume that Conjecture~\ref{conj:BAW} holds. Let $B$ be a block with abelian
 defect groups and $b$ its Brauer correspondent. Then $|\Irr(B)|=|\Irr(b)|$ and
 $|\IBr(B)|=|\IBr(b)|$.
\end{theorem}

Kn\"orr and Robinson have given a reformulation of the weight conjecture in
terms of certain simplicial complexes related to the $p$-local structure of
the group $G$. For this, let $\cP(G)$ denote
the set of chains $1<P_1<\ldots<P_l$ of $p$-subgroups of $G$. This induces a
structure of a simplicial complex on the set of non-trivial $p$-subgroups of
$G$. For $C=(1<P_1<\ldots<P_l)$ such a chain set $|C|=l$, the length of $C$,
and for
$B$ a $p$-block of $G$ let $B_C$ denote the union of all blocks $b$ of the
normaliser $\bN_G(C)$ with $b^G=B$. With this notation, Kn\"orr and Robinson
\cite[Thm.~3.8]{KR89} obtain the following reformulation:

\begin{theorem}[Kn\"orr--Robinson (1989)]   \label{thm:KR}
 The following two assertions are equivalent for a prime $p$:
 \begin{enumerate}
  \item[\rm(i)] Conjecture~\ref{conj:BAW} holds for all $p$-blocks of all
   finite groups;
  \item[\rm(ii)] for all $p$-blocks $B$ of all finite groups $G$ we have
   $$\sum_{C\in\cP(G)/\sim}(-1)^{|C|}|\IBr(B_C)| =0,$$
   where the sum runs over the chains in $\cP(G)$ up to $G$-conjugacy.
 \end{enumerate}
\end{theorem}

Here, in fact the set $\cP$ can also be replaced by the homotopy equivalent
sets of all chains of elementary abelian $p$-subgroups, or of all radical
$p$-subgroups, or by the set of chains in which all members are normal in the
larger ones.

By using M\"obius inversion it is possible from Theorem~\ref{thm:KR} to
describe the number of $p$-defect zero characters of $G$ in terms of local
subgroup information.

In the case of abelian defect groups, there is a strong relation between
Alperin's Weight Conjecture and the two previously introduced conjectures:

\begin{theorem}[Kn\"orr--Robinson (1989)]
 The following two assertions are equivalent for a prime $p$:
 \begin{enumerate}
  \item[\rm(i)] the Alperin--McKay Conjecture~\ref{conj:AM} holds for every
   $p$-block with abelian defect;
  \item[\rm(ii)] Alperin's Weight Conjecture~\ref{conj:BAW} holds for every
  $p$-block with abelian defect.
 \end{enumerate}
\end{theorem}

In fact, Kn\"orr--Robinson \cite[Prop.~5.6]{KR89} had to assume the
validity of the ``if''-direction of Conjecture~\ref{conj:H0} which is now
Theorem~\ref{thm:H0}.

The Alperin Weight Conjecture~\ref{conj:BAW} was proved by Isaacs and Navarro
\cite{IN95}
for $p$-solvable groups. It holds for all blocks with cyclic or non-abelian
metacyclic defect group by work of Brauer, Dade, Olsson and Sambale, see the
lecture notes \cite{Sam14}. It was shown to hold for groups of Lie type
in defining characteristic by Cabanes \cite{Ca88}, for $\fS_n$
and for $\GL_n(q)$ by Alperin and Fong \cite{AF90}, and by J.~An for certain
groups of classical type, see \cite{An94} and the references
therein. The latter proofs rely on an explicit determination of all weights
in the groups under consideration.

\subsection{Reductions}
As in the case of the Alperin--McKay conjecture, the Alperin weight
conjecture was first reduced in a non-block-wise form to some stronger
inductive statement (AWC) about finite simple groups by Navarro and Tiep
\cite{NT11} in 2011. In the same paper, they verified their inductive AWC
condition for example for groups of Lie type in their defining characteristic,
as well as for all simple groups with abelian Sylow 2-subgroups, while An and
Dietrich \cite{AD12} show it for sporadic groups. This reduction was then
refined by Sp\"ath \cite{Sp13b} to treat the block-wise version:

\begin{theorem}[Sp\"ath (2013)]
 The Alperin Weight Conjecture~\ref{conj:BAW} holds for all finite groups at
 the prime $p$ if all finite \emph{non-abelian simple} groups satisfy the
 so-called \emph{inductive block-wise Alperin weight condition (BAW)} at $p$.
\end{theorem}

Puig \cite{Pu11,Pu12} has announced another reduction of
Conjecture~\ref{conj:BAW} to nearly simple groups.

As in the case of the other inductive conditions, the \emph{inductive BAW
condition} for a simple group $S$ requires the existence of suitable equivariant
bijections at the level of the universal $p'$-covering group $G$ of $S$, this
time between $\IBr(B)$ and the weights attached to the block $B$ of $G$,
see \cite[Def.~4.1]{Sp13b} and also \cite{Sp16}.
In the same paper Sp\"ath shows that her inductive BAW condition holds for
various classes of simple groups, including the groups of Lie type in their
defining characteristic and for all simple groups with abelian Sylow
2-subgroup.

The inductive BAW condition has meanwhile been established by Breuer
\cite{BrWeb} for most sporadic simple groups, by the author \cite{Ma14} for
alternating groups, the Suzuki and the Ree groups, and by Schulte
\cite{Sch15} for the families of exceptional groups $G_2(q)$ and $\tw3D_4(q)$.
Koshitani and Sp\"ath \cite{KS15a} show that it holds for all blocks with
cyclic defect groups when $p$ is odd. 

For blocks $B$ with abelian defect groups Cabanes--Sp\"ath
\cite[Thm.~7.4]{CS13} and the author \cite[Thm.~3.8]{Ma14} have observed a
strong relation between the inductive BAW condition and the inductive
Alperin--McKay condition; we give here an even more general version from
\cite[Thm.~1.2]{KS15a}:

\begin{theorem}[Koshitani--Sp\"ath (2015)]   \label{thm:KS}
 Let $S$ be non-abelian simple with universal covering group $G$, $B$ a
 $p$-block of $G$ with abelian defect group $D$ and Brauer correspondent $b$
 in $\bN_G(D)$. Assume that the following hold:
 \begin{enumerate}
  \item[\rm(1)] The inductive Alperin--McKay condition holds for $B$ with
   respect to $M:=N_G(D)$ with a bijection $\Omega:\Irr_0(B)\rightarrow
   \Irr_0(b)$; and
  \item[\rm(2)] the decomposition matrix associated to
   $\Omega^{-1}(\{\chi\in\Irr(b)\mid D\le\ker(\chi)\})$ is lower uni-triangular
   with respect to some ordering of the characters.
 \end{enumerate}
 Then the inductive BAW condition holds for $B$ (considered as a block of
 $G/O_p(G))$.
\end{theorem}

This result highlights the importance of the existence of basic sets. Recall
that $X\subseteq\Irr(B)$ is a \emph{basic set for $B$} if the restrictions
to $p'$-elements of the $\chi\in X$ are linearly independent and span the
lattice $\ZZ\IBr(B)$ of Brauer characters. Such basic sets are known to exist
for groups of Lie type $\bG^F$ whenever the prime $p$ is good and does not
divide $|Z(\bG^F)|$: by a result of Geck and Hiss \cite{GH91}, $\cE(\bG^F,p')$
is a basic set for $\bG^F$, which moreover by definition is $\Aut(G)$-invariant.
Denoncin \cite{De15} has recently constructed such basic sets for the special
linear and unitary groups for all non-defining primes building on work
of Geck \cite{Ge91}.
It is an open question, formulated in \cite[(1.6)]{GH91} whether basic sets
exist for the blocks of all finite groups.

\begin{problem}
Construct natural $\Aut(G)_B$-invariant basic sets for blocks $B$ of finite
groups of Lie type $G$.
\end{problem}

Given an $\Aut(G)$-invariant basic set, condition~(2) of Theorem~\ref{thm:KS}
would be satisfied for example if the $p$-modular decomposition matrix of $G$
is lower unitriangular with respect to this basic set. This property is widely
believed to hold for groups of Lie type, and has been shown in a number
of important situations, for example by Gruber and Hiss \cite{GH97} if $G$ is
of classical type and the prime $p$ is \emph{linear} for $G$, as well as for
$\SL_n(q)$ and $\SU_n(q)$ by 
Kleshchev and Tiep \cite{KT09} and Denoncin \cite{De15}, respectively.

\begin{problem}
Show that decomposition matrices of finite reductive groups in non-defining
characteristic have uni-triangular shape.
\end{problem}

In the case of arbitrary defect it then still remains to determine the weights.
The weights of certain classical groups as well as of several series of
exceptional groups of Lie type of small rank have been determined by An and
collaborators, see e.g.~\cite{An94,ADH14}, but this has not resulted in a
general, type independent approach.

\begin{problem}
Give a generic description of weights of finite reductive groups, possibly
in the spirit of $d$-Harish-Chandra theory. Is there an analogue of Jordan
decomposition for weights?
\end{problem}

\section{Dade's Conjecture}   \label{sec:Dade}

Dade's Conjecture \cite{D92} extends the Kn\"orr--Robinson formulation in
Theorem~\ref{thm:KR} of the Alperin Weight Conjecture, and suggests a way to
count the characters of any defect in terms of the local subgroup structure.
It thus generalises both the McKay Conjecture~\ref{conj:AM} and the Alperin
Weight Conjecture~\ref{conj:BAW}. For this let us write
$$\Irr_d(B):=\{\chi\in\Irr(B)\mid \Ht(\chi)=d\}$$
for the irreducible characters in a block $B$ of height~$d$. Recall the set
$\cP(G)$ of chains of $p$-subgroups of $G$ from the previous section. The
so-called \emph{projective form} of Dade's conjecture claims:

\begin{conjecture}[Dade (1992)]   \label{conj:Dade}
 Let $B$ be a $p$-block of a finite group $G$. Then
 $$\sum_{C\in\cP(G)/\sim}(-1)^{|C|}|\Irr_d(B_C|\nu)| =0\qquad
   \text{for every }\nu\in\Irr(O_p(G))\text{ and }d\ge0,$$
 where the sum runs over chains in $\cP(G)$ up to $G$-conjugacy.
\end{conjecture}

As for the Kn\"orr--Robinson formulation of Alperin's Weight Conjecture, the
set $\cP$ of chains may be replaced by chaines involving only elementary
abelian $p$-subgroups, or only radical $p$-subgroups.
 
Dade's Conjecture was proved for $p$-solvable groups by Robinson \cite{Ro00}.
An has shown Dade's conjecture for general linear and unitary groups in
non-defining characteristic, and for various further groups of Lie type of
small rank, see e.g.~\cite{An01}.
Recently, in a tour de force Sp\"ath \cite{Sp15} managed to reduce a suitable
form of Dade's conjecture to a statement on simple groups:

\begin{theorem}[Sp\"ath (2015)]
 Dade's projective Conjecture~\ref{conj:BAW} holds for all finite groups at
 the prime $p$, if all finite \emph{non-abelian simple} groups satisfy the
 so-called \emph{character triple conjecture} at $p$.
\end{theorem}

The character triple conjecture (see \cite[Conj.~1.2]{Sp15}) is a statement
about chains in $\cP$ similar to Dade's projective conjecture, but as in the
previous inductive conditions it also involves the covering groups and the
action of automorphisms. It has been proved for blocks with cyclic defect,
the blocks of sporadic quasi-simple groups except for the baby monster $B$
and the monster $M$ at $p=2$, and for  $\PSL_2(q)$ \cite[Thm.~9.2]{Sp15}.

\begin{problem}
Find a generic way to describe $p$-chains in finite reductive groups.
\end{problem}

\section{Brou\'e's Abelian Defect Group Conjecture}
An attempt to give a structural explanation for all of the ``numerical''
local-global conjectures mentioned so far is made by Brou\'e's conjecture at
least in the case of blocks with abelian defect groups. Recall that the
Alperin--McKay Conjecture~\ref{conj:AM} relates character degrees of a
$p$-block $B$ of a finite group $G$ with defect group $D$ to those of a
Brauer corresponding block $b$ of $\bN_G(D)$. Brou\'e \cite{Br90} realised
that this numerical relation would be a consequence of a (conjectural) intimate
structural relation between the module categories of the $\cO$-algebras $B$
and $b$:

\begin{conjecture}[Brou\'e (1990)]   \label{conj:Br}
 Let $B$ be a block of a finite group with defect group $D$ and $b$ its Brauer
 corresponding block of $\bN_G(D)$. Then the bounded derived module categories
 of $B$ and of $b$ are equivalent.
\end{conjecture}

Brou\'e shows that the validity of his conjecture for a block $B$ would imply
the Alperin--McKay conjecture as well as the Alperin weight conjecture for $B$.

Brou\'e's conjecture holds for all blocks of $p$-solvable groups, since in
this case by Dade \cite{D80} and Harris--Linckelmann \cite{HL00}, the two
blocks in question are in fact Morita equivalent. It also holds for blocks with
cyclic or Klein four defect groups by work of Rickard \cite{Ri89,Ri96}. Using
the classification of finite simple groups it has been shown for principal
blocks with defect group $C_3\times C_3$ by Koshitani and Kunugi \cite{KK02},
and it has been verified for many blocks with abelian defect of sporadic
simple groups.

In their landmark paper \cite{CR1} Chuang and Rouquier have given a proof of
Brou\'e's conjecture for $\fS_n$ and for $\GL_n(q)$, building on
previous results of Chuang and Kessar \cite{CK02} and Turner \cite{Tu02} who
obtained derived equivalences for a very particular class of blocks, the
so-called Rouquier blocks. Dudas, Varagnolo and Vasserot \cite{DVV} have
constructed derived equivalences between blocks of various finite unitary
groups which together with a result of Livesey \cite{Liv15} provides a
verification of Brou\'e's conjecture for $\GU_n(q)$ for linear primes.

Brou\'e and the author \cite{BMM} proposed a more precise form of
Conjecture~\ref{conj:Br} in the case of unipotent blocks of finite groups of
Lie type in terms of the $\ell$-adic cohomology of Deligne--Lusztig varieties.
This version has recently been explored by Craven and Rouquier \cite{CR2}
using the concept of perverse equivalences.

Despite of these partial successes, in contrast to the situation for the other
conjectures stated earlier, there is no general reduction theorem for Brou\'e's
conjecture to a condition on simple groups. A further challenge lies in the
fact that currently Brou\'e's conjecture is only formulated for blocks with
abelian defect groups, and it remains unclear how a generalisation to blocks
of arbitrary defect might look like.


\end{document}